\documentclass[11pt,a4paper]{article}
\usepackage[utf8]{inputenc}
\usepackage[english]{babel}
\usepackage[T1]{fontenc}

\usepackage{amssymb}
\usepackage{amsmath}
\usepackage{amsthm}
\usepackage{framed}
\usepackage{braket}
\usepackage{bm}
\usepackage{mathrsfs}
\usepackage{accents}
\usepackage{tocloft}
\usepackage{graphicx}
\usepackage{tikz}
\usepackage{url}
\usepackage{color}
\usepackage{xifthen}
\usepackage{xcolor}
\usepackage{framed}
\usepackage{mathtools}
\usepackage[explicit]{titlesec}
\usepackage{mdframed}
\usepackage{geometry}
\usepackage{enumerate}
\usepackage{hyperref}

\usetikzlibrary{positioning}
\usetikzlibrary{calc}
\usetikzlibrary{decorations.pathreplacing}
\usetikzlibrary{decorations.fractals}
\usetikzlibrary{cd}

\newcommand{\scrN}{\mathcal{N}}

\newcommand{\R}{\mathbb{R}}

\newcommand{\non}{\operatorname{non}}
\newcommand{\cov}{\operatorname{cov}}
\newcommand{\add}{\operatorname{add}}
\newcommand{\cof}{\operatorname{cof}}
\newcommand{\nul}{\mathcal{N}}
\newcommand{\meager}{\mathcal{M}}

\newcommand{\diam}{\mathrm{diam}}

\newcommand{\Haus}{\mathcal{H}}
\newcommand{\HDZ}{\mathsf{HDZ}}

\newcommand{\I}{\mathcal{I}}

\newcommand{\seq}[1]{{\langle#1\rangle}}
\DeclarePairedDelimiter\abs{\lvert}{\rvert}

\renewcommand\subset{\subseteq}

\renewcommand{\le}{\leqslant}
\renewcommand{\ge}{\geqslant}

\newcommand{\needtocheck}[1][]{%
	\ifthenelse{\equal{#1}{}}{%
		\textcolor{blue}{[NeedToCheck]}%
	}{%
		\textcolor{blue}{[NeedToCheck: #1]}%
	}%
}

\newcommand{\todo}[1][]{%
	\ifthenelse{\equal{#1}{}}{%
		\textcolor{red}{[TODO]}%
	}{%
		\textcolor{red}{[TODO: #1]}%
	}%
}

\theoremstyle{definition}
\newtheorem{thm}{Theorem}[section]
\newtheorem*{thm*}{Theorem}

\newtheorem*{defi*}{Definition}
\newtheorem{lem}[thm]{Lemma}
\newtheorem*{lem*}{Lemma}

\newtheorem*{fact*}{Fact}

\newtheorem*{prop*}{Proposition}

\newtheorem*{exm*}{Example}

\newtheorem*{rmk*}{Remark}

\newtheorem*{cor*}{Corollary}
\newtheorem*{notation*}{Notation}

\usepackage{xpatch}
\newcounter{proofdepth}

\xpretocmd{\proof}{%
	\stepcounter{proofdepth}%
	\ifnum\value{proofdepth}=1
	\else
	\fi
}{}{}

\xapptocmd{\endproof}{%
	\addtocounter{proofdepth}{-1}%
}{}{}
\newcounter{claim}[thm]

\newtheorem{innerclaim}{Claim}

\definecolor{reasontext}{rgb}{0,0,0}
\definecolor{reasonbg}{rgb}{0.9,0.9,0.9}

\usepackage{csquotes}
\usepackage[backend=biber,style=alphabetic,sorting=nty,doi=false,isbn=false,url=false,eprint=false]{biblatex}
\renewbibmacro{in:}{}
\definecolor{myred}{RGB}{230, 150, 100}
\definecolor{mygreen}{RGB}{150, 230, 100}
\definecolor{myblue}{RGB}{130, 220, 220}
\hypersetup{linkbordercolor=myred,citebordercolor=mygreen,urlbordercolor=myblue}

\title{The cardinal invariants of the $s$-dimensional Hausdorff null ideal are independent of $s$}
\author{Tatsuya Goto}
\date{\today}

\begin{document}
\maketitle

\begin{abstract}
For $0 < s < 1$, let $\scrN^s$ denote the ideal of sets of $s$-dimensional Hausdorff measure zero.
We prove that $\cov(\scrN^s) = \cov(\scrN^t)$ and $\non(\scrN^s) = \non(\scrN^t)$ for all $0 < s, t < 1$.
The proof uses Ren and Wang's solution of the Furstenberg set conjecture.
This answers open questions of Shelah--Stepr\=ans and Elekes--Stepr\=ans.
\end{abstract}

\section{Introduction}

The four cardinal invariants $\add$, $\cov$, $\non$, and $\cof$ of the meager ideal $\meager$, the Lebesgue null ideal $\nul$, and the ideal $K_\sigma$ generated by the compact subsets of Baire space have been extensively studied~\cite{Blass2010}.
There is also growing interest in the study of ideals other than these three.

Hausdorff measure and dimension quantify the size of Lebesgue null sets~\cite{Mattila1995}.
For each $0 < s < 1$, one can consider the $s$-dimensional Hausdorff null ideal $\scrN^s$.
The following results are known.
\begin{enumerate}
    \item $\add(\scrN^s) = \add(\nul)$ for every $0 < s < 1$.
    The same holds for $\cof$ (Fremlin \cite[534B]{Fremlin2008}).
    \item It is consistent that $\non(\scrN^s) < \non(\nul)$ for all $0 < s < 1$ (Shelah--Stepr\=ans \cite{ShelahSteprans2005}, with a simplified proof by Goto \cite{Goto2025}).
    \item It is consistent that $\cov(\nul) < \cov(\scrN^s)$ for all $0 < s < 1$ (Elekes--Stepr\=ans \cite[Corollary~3.3]{ElekesSteprans2019}).
\end{enumerate}

Shelah--Stepr\=ans asked whether it is consistent that $\non(\scrN^s) < \non(\scrN^t)$ for some $0 < s < t < 1$~\cite[Question~7.1]{ShelahSteprans2005}.
Similarly, Elekes--Stepr\=ans asked whether it is consistent that $\cov(\scrN^t) < \cov(\scrN^s)$ for some $0 < s < t < 1$~\cite[Problem~5.4]{ElekesSteprans2019}.
In this paper, we answer both questions negatively.

Let $X$ be a metric space, $A \subset X$, and $s,\delta > 0$.
The $\delta$-approximation to the $s$-dimensional Hausdorff measure of $A$ is defined by
\[
\Haus^s_\delta(A)
= \inf\left\{
    \sum_{n=0}^\infty \diam(C_n)^s :
    \begin{gathered}
        C_n \subseteq X \quad (n \in \omega),\\
        A \subseteq \bigcup_{n \in \omega} C_n,\quad
        (\forall n)\bigl(\diam(C_n) \le \delta\bigr)
    \end{gathered}
\right\}.
\]
The $s$-dimensional Hausdorff measure of $A$ is defined by
\[
\Haus^s(A) = \lim_{\delta \to 0} \Haus^s_\delta(A).
\]
The Hausdorff dimension of $A$ is defined by
\[
\dim_\mathrm{H}(A)
= \sup \{ s : \Haus^s(A) = \infty \}
= \inf \{ s : \Haus^s(A) = 0 \}.
\]
The last equality is well known and easy to prove.

We equip Cantor space $2^\omega$ with the metric
\[
d(x,y) =
\begin{cases}
    0 & \text{if $x=y$,}\\
    2^{-\min\{n:x(n)\ne y(n)\}} & \text{otherwise.}
\end{cases}
\]

For $0 < s < 1$, let $\scrN^s$ be the ideal of subsets of Cantor space $2^\omega$ of $s$-dimensional Hausdorff measure zero.
Let $\HDZ$ be the ideal of subsets of $2^\omega$ of Hausdorff dimension zero.

Our main tool is K.~Ren and H.~Wang's solution of the Furstenberg set conjecture~\cite[Theorem~1.1]{RenWang2023}.

Let $0 < s \le 1$ and $0 < t \le 2$.
A set $E \subset \R^2$ is an $(s,t)$-Furstenberg set if there is a subset $L$ of the space $\mathcal{A}(2,1)$ of affine lines in the Euclidean plane such that
$\dim_\mathrm{H}^{\mathcal{A}(2,1)}(L) \ge t$ and
$\dim_\mathrm{H}(E \cap \ell) \ge s$ for every $\ell \in L$.
Ren and Wang's theorem states that
\[
\dim_\mathrm{H}(E)
\ge \min\left\{s+t,\frac{3s+t}{2},s+1\right\}.
\]
Intuitively, a planar set that has large intersections with many lines must itself be large.

Here we equip $\mathcal{A}(2,1)$ with the metric
\[
d(\ell,\ell') = \sin\theta + |p-p'|,
\]
as in Section~3.16 of \cite{Mattila1995}.
Here $\theta\in[0,\pi/2]$ is the smaller angle between the two lines, and $p,p'$ are the points on the respective lines closest to the origin.
The map
\begin{equation}\label{eq:mapT}
T:[0,1]^2\to\mathcal{A}(2,1),\qquad T(u,v)=\ell_{u,v},
\end{equation}
is bi-Lipschitz~\cite[Section~6.1]{DabrowskiOrponenVilla2022}, where
\[
\ell_{u,v} = \left\{(x,y)\in\R^2:y=\frac{u}{2}x+\frac{v}{2}\right\}.
\]

Using dimension estimates for Furstenberg sets, Y.~Wang proved a zero--one law: every Turing ideal $\mathcal{I}\subseteq 2^\omega$ satisfies $\dim_\mathrm{H}\mathcal{I}\in\{0,1\}$~\cite[Corollary~1.2]{Wang2026}.
The assertion about $\non$ in our main result also follows directly from Wang's theorem.

Finally, we recall the definitions of the four cardinal invariants of an ideal.
For an ideal $\I$ on a set $X$:
\begin{enumerate}
    \item $\add(\I)$, the additivity of $\I$, is the smallest cardinality of a family $F\subseteq\I$ such that $\bigcup F\notin\I$.
    \item $\cov(\I)$, the covering number of $\I$, is the smallest cardinality of a family $F\subseteq\I$ such that $\bigcup F=X$.
    \item $\non(\I)$, the uniformity of $\I$, is the smallest cardinality of a set $A\subseteq X$ such that $A\notin\I$.
    \item $\cof(\I)$, the cofinality of $\I$, is the smallest cardinality of a family $F\subseteq\I$ such that for every $A\in\I$ there is a $B\in F$ with $A\subset B$.
\end{enumerate}

\section{\texorpdfstring{$\non(\scrN^s) = \non(\scrN^t)$}{non(Ns) = non(Nt)}}

\begin{thm}\label{thm:non}
Let $0 < s < 1$. Then $\non(\scrN^s) = \non(\HDZ)$.
\end{thm}
\begin{proof}
Since $\HDZ \subset \scrN^s$, the inequality $\non(\HDZ) \le \non(\scrN^s)$ is immediate.
We prove the reverse inequality.
Choose $A\notin\HDZ$ with $\abs{A}=\non(\HDZ)$, and let $I(A)$ be the Turing ideal generated by $A$.
For each real $x$, there are at most countably many reals Turing reducible to $x$.
Thus $\abs{I(A)}=\abs{A}=\non(\HDZ)$.
By Wang's theorem, $I(A)$ has Hausdorff dimension $1$.
Hence $I(A)\notin\scrN^s$, and therefore
\[
\non(\scrN^s) \le \abs{I(A)} = \non(\HDZ).
\]
\end{proof}

\section{\texorpdfstring{$\cov(\scrN^s) = \cov(\scrN^t)$}{cov(Ns) = cov(Nt)}}

Define $\Phi:\R\to\R$ by
\[
\Phi(s) = \min\left\{\frac{5}{4}s,\frac{1+s}{2}\right\}.
\]
For $E\subset\R^2$ and $0 < s \le 1$, put
\[
P_s(E) = \{(u,v)\in[0,1]^2:\dim_\mathrm{H}(E\cap\ell_{u,v})\ge s\}.
\]

\begin{lem}\label{lem:cor_of_ren_wang}
Let $E\subset\R^2$ and $0 < s \le 1$.
If $\dim_\mathrm{H}(P_s(E))\ge 2s$, then $\dim_\mathrm{H}(E)\ge 2\Phi(s)$.
\end{lem}
\begin{proof}
Using the map $T$ in \eqref{eq:mapT}, put $L=T[P_s(E)]$.
Then $L$ witnesses that $E$ is an $(s,2s)$-Furstenberg set.
By the Ren--Wang theorem,
\[
\dim_\mathrm{H}(E)
\ge \min\left\{3s,\frac{3s+2s}{2},s+1\right\}
= \min\left\{\frac{5s}{2},s+1\right\}
= 2\Phi(s).
\]
\end{proof}

\begin{thm}
For all $0 < s,t < 1$, we have $\cov(\scrN^s)=\cov(\scrN^t)$.
\end{thm}
\begin{proof}
It suffices to consider the case $0 < s < t < \Phi(s)$.
Indeed, once the result is established in this case, for arbitrary $0 < s < t < 1$ we can choose finitely many points
$s=s_0<\cdots<s_m=t$ with $s_{i+1}<\Phi(s_i)$.
Applying the equality to each adjacent pair yields the equality for $s$ and $t$.
Here we use the fact that $\Phi(s)-s>0$ for every $0 < s < 1$.

Suppose that $\kappa=\cov(\scrN^t)<\cov(\scrN^s)$.
Choose a sequence $\seq{A_\alpha:\alpha<\kappa}$ such that
$A_\alpha\in\scrN^t$ and $\bigcup_{\alpha<\kappa}A_\alpha=2^\omega$.

To pass from Cantor space to the Euclidean plane, define $g:2^\omega\to[0,1]^2$ by
\[
g(x) = \left(\sum_{n\in\omega}\frac{x(2n)}{2^{n+1}},
             \sum_{n\in\omega}\frac{x(2n+1)}{2^{n+1}}\right).
\]
By considering the H\"older properties of this map, we obtain
\[
\begin{aligned}
\dim_\mathrm{H}g[A] &= 2\dim_\mathrm{H}A
    &&(A\subset 2^\omega),\\
\dim_\mathrm{H}g^{-1}[B] &= \tfrac12\dim_\mathrm{H}B
    &&(B\subset[0,1]^2).
\end{aligned}
\]

For each $\alpha<\kappa$, put $E_\alpha=g[A_\alpha]$ and
\[
P_\alpha = \{(u,v)\in[0,1]^2:\dim_\mathrm{H}(E_\alpha\cap\ell_{u,v})\ge s\}.
\]
We have $\dim_\mathrm{H}(E_\alpha)\le 2t<2\Phi(s)$.
By the contrapositive of Lemma~\ref{lem:cor_of_ren_wang},
$\dim_\mathrm{H}(P_\alpha)<2s$.
Thus $\Haus^s(g^{-1}[P_\alpha])=0$.
Since $\kappa<\cov(\scrN^s)$, it follows that
\[
\bigcup_{\alpha<\kappa}g^{-1}[P_\alpha]\ne 2^\omega.
\]
Hence there is an $x^*\in 2^\omega$ such that $g(x^*)\notin P_\alpha$ for every $\alpha<\kappa$.
Write $g(x^*)=(u^*,v^*)$.

Put
\[
B_\alpha = \left\{x\in[0,1]:
    \left(x,\frac{u^*}{2}x+\frac{v^*}{2}\right)\in E_\alpha\right\}.
\]
Since $(u^*,v^*)\notin P_\alpha$ and the projection from $\ell_{u^*,v^*}$ onto the $x$-axis is Lipschitz, we have $\dim_\mathrm{H}(B_\alpha)<s$.
On the other hand, the sets $E_\alpha$ ($\alpha<\kappa$) cover $[0,1]^2$, so the sets $B_\alpha$ ($\alpha<\kappa$) cover $[0,1]$.
To obtain a contradiction to $\kappa<\cov(\scrN^s)$, we pass back to Cantor space.
Define $h:2^\omega\to[0,1]$ by
\[
h(x) = \sum_{n\in\omega}\frac{x(n)}{2^{n+1}}.
\]
Again, considering the H\"older properties gives
\[
\begin{aligned}
\dim_\mathrm{H}h[A] &= \dim_\mathrm{H}A
    &&(A\subset 2^\omega),\\
\dim_\mathrm{H}h^{-1}[B] &= \dim_\mathrm{H}B
    &&(B\subset[0,1]).
\end{aligned}
\]
Now $\bigcup_{\alpha<\kappa}h^{-1}[B_\alpha]=2^\omega$, while
$\dim_\mathrm{H}(h^{-1}[B_\alpha])<s$ for each $\alpha<\kappa$.
This contradicts $\kappa<\cov(\scrN^s)$ and completes the proof.
\end{proof}

\section{\texorpdfstring{$\non(\scrN^s) = \non(\scrN^t)$}{non(Ns) = non(Nt)} revisited}

Theorem~\ref{thm:non} can also be proved using the Ren--Wang theorem.

\begin{thm}
For all $0 < s,t < 1$, we have $\non(\scrN^s)=\non(\scrN^t)$.
\end{thm}
\begin{proof}
As before, it suffices to consider the case $0 < s < t < \Phi(s)$.
Choose $X\subset 2^\omega$ such that $X\notin\scrN^s$ and $\abs{X}=\non(\scrN^s)$.
Using the maps $g:2^\omega\to[0,1]^2$ and $h:2^\omega\to[0,1]$ from the preceding proof, put
\[
E = \left\{\left(x,\frac{u}{2}x+\frac{v}{2}\right):
    x\in h[X],\ (u,v)\in g[X]\right\}.
\]
The family of lines parametrized by $g[X]$ has Hausdorff dimension at least $2s$, and each line contains a bi-Lipschitz image of $h[X]$.
By Lemma~\ref{lem:cor_of_ren_wang}, $\dim_\mathrm{H}(E)\ge 2\Phi(s)$.
Put $Y=g^{-1}[E]$.
Then $\dim_\mathrm{H}(Y)\ge\Phi(s)>t$, and $\abs{Y}=\abs{X}$. Although $g$ is not injective, it is at most four-to-one.
Consequently,
\[
\non(\scrN^t) \le |Y| = |X| = \non(\scrN^s).
\]
\end{proof}

\section*{Declaration of generative AI}
During the preparation of this work, the author used OpenAI ChatGPT and Codex with
the \texttt{gpt-6.2-astra} and  \texttt{gpt-6.2-sol} models to assist with
literature searches, the discovery and development of proofs, language
editing, and \LaTeX{} drafting. The author independently checked,
reviewed, and revised all AI-assisted material and takes full
responsibility for the content of the paper.

\printbibliography

\end{document}